\documentclass[11pt]{article}
\usepackage[pctex32]{graphics}
\usepackage{amssymb}
\usepackage{latexsym}
\usepackage{epsfig} 
\usepackage{multirow} 
\usepackage{pstcol}
\usepackage{float}
\usepackage{epstopdf}
\usepackage{mathrsfs}
\usepackage{soul}

\newcommand\bs{\boldsymbol} 

\usepackage{amsmath}
\usepackage{amssymb}
\usepackage{relsize}
\usepackage{multirow}
\usepackage{algorithm}
\usepackage{algpseudocode}
\usepackage{graphicx}
\usepackage{multicol} 

\usepackage{times}
\usepackage{epstopdf}

\definecolor{Red}{rgb}{1,0.,0.}
\usepackage{fancyhdr}
\usepackage{setspace}
\usepackage{hyperref}
\usepackage{mathtools}
\usepackage[dvipsnames]{xcolor}

\usepackage{placeins}
\usepackage{amsopn}
\newcommand{\argmin}{{\rm arg}\min}

\title{Sparse Dictionary-Based Solution of Dynamic Inverse Problems}

\author{Aidan Mason-Mackay$^{1,3}$ \and Daniela Calvetti$^2$ \and Erkki Somersalo$^2$ \and Antti Aarnio$^3$ \and Mikko Kettunen$^3$ \and Ekaterina Paasonen$^3$ \and Olli Gröhn$^3$ 
\and Ville Kolehmainen$^1$}

\date{$^1$Department of Technical Physics, Faculty of Science, Forestry and Technology, University of Eastern Finland, FI-70211 Kuopio, Finland  \\
$^2$Department of Mathematics, Applied Mathematics and Statistics, Case Western Reserve University, Cleveland, OH 4410, USA  \\
$^3$A.I Virtanen Institute for Molecular Sciences, Faculty of Health Sciences, University of Eastern Finland, FI-70211 Kuopio, Finland
}

\begin{document}

\maketitle

% REQUIRED
\begin{abstract}
In ill-posed dynamic inverse problems expected spatial features and temporal correlation between frames can be leveraged to improve the quality of the computed solution, in particular when the available data are limited and the dimensionality of the unknown is large. One way to take advantage of the spatial and temporal traits believed to characterize the solution is to encode them into the entries of a dictionary, and to seek the solution as a sparse linear combination of the dictionary atoms. To promote a vector of coefficients with mostly vanishing entries, we consider a stochastic extension of the dictionary coding problem model with a random hierarchical sparsity promoting prior. We compute the Maximum A Posteriori (MAP) estimate of the coefficient vector using the Iterative Alternating Sequential Algorithm (IAS), which has been demonstrated to efficiently solve inverse problems with minimal need for parameter tuning. The proposed methodology is tested on real-world dynamic Computed Tomography and MRI datasets, where it is compared to the popular Alternating Direction Method of Minimizers (ADMM). The computed examples show the that proposed methodology is competitive with the ADMM for compressed sensing, with a significantly lower sensitivity to hyper-parameter selection. 
\end{abstract}

% REQUIRED
%\begin{keywords}
%Dynamic Inverse Problems, Sparsity Promotion, MRI, Computed Tomography, Hierarchical Bayesian Methods
%\end{keywords}

% REQUIRED
%\begin{MSCcodes}
%47A52, 65J22, 94A08
%\end{MSCcodes}

\section{Introduction}
In dynamic inverse problems, the goal is to estimate a time series of parameters from indirect measurements. Such problems arise in many settings, such as medical imaging, including functional Magnetic Resonance Imaging (fMRI), Computed Tomography (CT), Single Photon Computed Tomography (SPECT) and Positron Emission Tomography (PET).  The time-dependent parameters in these problems are typically spatially distributed, and consequently, the  models are often high-dimensional, and computational efficiency is an important consideration. Moreover, the maximum achievable sampling rate is often limited by physical constraints, and, as such, one may need to acquire heavily under-sampled data in order to achieve the desired temporal resolution. Due to the underdetermined nature of the problem, as well as the presence of measurement noise, dynamic inverse problems of this type are characteristically  ill-posed, and a priori information must be incorporated into the reconstruction process to arrive at a meaningful estimate of the unknown parameters. 

There is a large literature on methods for solving ill-posed dynamic inverse problems, see, e.g., \cite{hauptmann_image_2020, schuster_dynamic_2018}. Many approaches exploit temporal correlations between frames, including  variational techniques with temporal regularization terms, such as Tikhonov regularization \cite{schmitt_efficient_2002}, total generalized variation \cite{bredies_total_2010}, or low-rank promotion \cite{gao_robust_2011}, which may be combined with spatial penalty terms. The inverse problem can be formulated in a classical deterministic or Bayesian framework \cite{calvetti_inverse_2018}, and there is also a  growing literature on data-driven reconstruction techniques \cite{arridge_solving_2019}. 

The solution of a dynamic ill-posed inverse problem from few data can be formulated as a sparse dictionary coding problem, where the unknown quantity of interest is represented as a linear combination of a few of an ensemble of solution templates, referred to as the atoms of a dictionary. 
In this paper we focus our attention on  dictionary-based solution methods for dynamic inverse problems where the dictionaries are designed to account simultaneously for  the spatial and temporal priors, and the solution is expressed in terms of only few dictionary atoms. More specifically, the atoms of the dictionary are stored as the columns of a matrix, and the prior information about the unknown is encoded in the dictionary entries via appropriate matrix operations. Typically, the dictionary reflects the nature of the underlying application: for example, in \cite{calvetti_dictionary_2025}, the dictionary is formed based on an anatomical atlas of the brain, while in \cite{bocchinfuso_bayesian_2025} the dictionary atoms are the computed responses by the governing equations corresponding to different input vectors.
The assumption that only a few dictionary items are needed to explain the unknown motivates the adoption of a Bayesian prior model promoting sparsity in the vector of coefficients defining the linear combination. 
In summary, the dictionary coding is formally a linear inverse problem, and its sparse solution in this work is computed using the The Iterative Alternating Sequential Algorithm (IAS), a sparsity-promoting reconstruction method, leveraging the formulation of the problem as a hierarchical Bayesian model which has been used recently in sparse dictionary coding problems; see e.g. \cite{pragliola2021overcomplete, bocchinfuso_bayesian_2025}.  

The IAS method, a two-phase block descent algorithm designed for the computation of the MAP estimate based on hierarchical prior models, has been shown to converge to the unique solution of its associated strictly convex objective function; see e.g., \cite{calvetti2023bayesian, calvetti2019hierachical}. IAS has been shown to be consistently computationally efficient and capable of producing high quality reconstructions from both synthetic and real-world data in a variety of application areas, including image reconstruction and denoising \cite{pragliola2021overcomplete}, optical flow velocimetry \cite{jassal_bayesian_2025}, cell membrane permeability estimation \cite{bocchinfuso_bayesian_2025},  PET imaging \cite{bardsley_hierarchical_2010}, electrical impedance tomography and optical tomography \cite{manninen2022sparsity,calvetti2025sparsity}, and brain activity mapping using both MEG and EEG data \cite{calvetti_ias-meeg_2023, calvetti_dictionary_2025}.
The convergence properties of the IAS algorithm are well understood \cite{calvetti2019hierachical,calvetti2020sparse}: For suitable choices of hyperpriors, the cost function associated with the hierarchical Bayesian formulation is strictly convex, and the IAS algorithm convergences quadratically to the zero background. A discussion of the physical interpretation of the IAS hyper-parameters and automatic hyper-parameter selection schemes can be found in \cite{calvetti_ias-meeg_2023}.

For problems where the signal is incoherent and sparse, as is often the case in applications such as MRI, compressed sensing (CS) based algorithms are very popular \cite{uecker2015berkeley, fessler_optimization_2019}, and routinely are used by commercial imaging machines \cite{forman2016compressed}. In the compressed sensing framework, perfect reconstructions from suitably under-sampled, noisy data are possible with high probability, provided that the parameter of interest is sufficiently sparse in the given transform domain \cite{candes_robust_2006}. Dictionary methods have been used in the CS framework, including techniques of learning the sparsifying dictionaries 
from data \cite{caballero2014dictionary,huang2014bayesian} as well as using dictionaries for quantitative imaging
\cite{ma2013magnetic}. CS based algorithms typically include one or more sparsifying regularizing penalty functionals, and it is well understood that regularization, needed to overcome the ill-posedness of compressed sensing problems, has a significant effect on the reconstruction quality, and requires empirical tuning of hyper-parameters \cite{forman2016compressed}. If the problem is over-regularized, there is the risk of dimming or removing the signal of interest, while if the problem is under-regularized, the recovered solutions may be dominated by amplified noise. Selection of compressed sensing hyper-parameters continues to be an active area of research, see for example \cite{wang_regularization-agnostic_2021}. Popular sparsifying transforms used in compressed sensing include temporal or spatial finite differences and wavelet transforms \cite{rani_systematic_2018}. The cost function associated with a compressed sensing problem may be non-smooth and convex, and is commonly minimised using the Alternating Direction Method of Minimizers (ADMM) \cite{boyd_distributed_2010}. In the ADMM method, a third hyper-parameter is generally tuned empirically following some heuristics. The ADMM method has been shown to be globally convergent in the context of CS, with sublinear to linear convergence rates \cite{boyd_distributed_2010}. The convergence rate of the ADMM method can be improved by adaptive penalty selection strategies such as residual balancing \cite{wohlberg_admm_2017}, at the cost of potentially weakening  convergence guarantees. The ADMM method has been shown empirically to produce solutions of modest accuracy in a few tens of iterations, while requiring much longer for high-accuracy solutions \cite{boyd_distributed_2010}.

In this paper, the proposed IAS-based approach is shown to be a competitive alternative for CS based algorithms both in computational efficiency and quality of the results. The algorithm
is compared to the CS ADMM approach on a dynamic X-ray Computed Tomography (CT) problem using experimental data from a phantom target, and on a Dynamic-Contrast Enhanced MRI (DCE-MRI) using experimental data from a small animal study of imaging brain glioma, focusing on reconstruction accuracy, hyper-parameter sensitivity, and runtimes. In particular, while the quality of the reconstructions computed by the two methods is similar under optimal parameter selection, the IAS-based approach is much less sensitive to  hyperparameter selection and requires less time to converge to the solution. 

The paper is organized as follows. Section \ref{sec:problemDefinition} defines the dynamic inverse problem, Sections \ref{subsec:dictionary}-\ref{sec:IAS} describe the proposed dictionary formulation and reconstruction strategy, results are presented in Section \ref{sec:results}, and Section \ref{sec:conclusions} contains concluding remarks.

\section{Problem Definition}
\label{sec:problemDefinition}
Consider a dynamic inverse problem discretized in space and time, in which the aim is to estimate a time series of unknown unobservable vectors based on a time series of indirect observations. Let 
\begin{equation}\label{TimeSeriesX}
X = [\bs{x}^{(1)} \; \bs{x}^{(2)} \;  \hdots \; \bs{x}^{(N_T)}] \in \mathbb{C}^{N_P \times N_T}
\end{equation}
be the matrix collecting in each of its columns $\bs{x}^{(t)}$ the state of the unknown vector at the time instance $t$, and 
\begin{equation}\label{observablesB} 
B = [\bs{b}^{(1)}  \; \bs{b}^{(2)} \;  \hdots \;  \bs{b}^{(N_T)}] \in \mathbb{C}^{N_D \times N_T},
\end{equation}
be the matrix where each column $\bs{b}^{(t)}$ contains the observed quantity at  time instance $t$.
It is assumed that at a given time instant $t$, $1\leq t\leq N_T$, the vectors $\bs{x}^{(t)}$ and $\bs{b}^{(t)}$ are related by the linear forward model
\begin{align}\label{forward}
    \bs{b}^{(t)} = F^{(t)} \bs{x}^{(t)} + \bs{\epsilon}^{(t)},
\end{align}
where the matrix $F^{(t)} \in \mathbb{C}^{N_D \times N_P}$ maps the unknown to data space, and $\bs{\epsilon}^{(t)} \in \mathbb{C}^{N_P}$ represents additive noise that is assumed here to be stationary, zero mean and Gaussian, i.e. 
$\bs{\epsilon}^{(t)} \sim{\mathcal N}(0,\Gamma)$, where $\Gamma\in \mathbb{C}^{N_T\times N_T}$ is Hermitian symmetric positive definite.  While it is convenient to formulate the forward problem in terms of the matrix $F^{(t)}$, in many practical applications such as MRI, CT, and PET,  the matrix is not explicitly formed, but  computations requiring multiplication of given vectors by $F^{(t)}$ or its adjoint can be  implemented in a matrix-free fashion.
In general, both the unknowns to be estimated and the observed quantities are assumed to be complex vectors, the restriction to real vectors being a special case.

Denoting by $vec$ the operator stacking the columns of a matrix from left to right into a vector, we define 
$\bs{x} = vec \left( X \right)$, $\bs{b}=vec\left( B \right)$, $\bs{\epsilon}= vec \left( [\bs{\epsilon}^{(1)} \; \bs{\epsilon}^{(2)} \; \dots \; \bs{\epsilon}^{(N_T)}] \right)$, and introduce the block diagonal matrix $F \in \mathbb{C}^{M  \times N}$,
\begin{align*}
F = \begin{bmatrix}
    F^{(1)} & & & \\
     & F^{(2)} &  &  \\
     & & \ddots &  \\
     & &  & F^{(N_T)}
    \end{bmatrix},
\end{align*}
where $M=(N_P N_T)$ and $N=(N_D N_T)$. 
With these notations, we can express the forward models (\ref{forward}) as the single linear system 
\begin{align} \label{eqn:ForwardModel}
    \bs{b} = F \bs{x} + \bs{\epsilon}.
\end{align}

We are now ready to reformulate the inverse problem of estimating $\bs{x}$, given observed data $\bs{b}$ and the model  (\ref{eqn:ForwardModel}) in probabilistic terms. Following the Bayesian paradigm, the unknown $\bs{x}$ is modeled as a random variable $\bs{X}$ with realizations denoted by $\bs{x}$. In line with the statement of the problem, we assume that the noise vector follows a Gaussian distribution ${\mathcal N}(0,I_{N_T}\otimes\Gamma)$. Consequently, it follows that if $\bs{X}=\bs{x}$ is given, the observation vector is also a a random variable $\bs{B}$ and its conditional probability $\pi_{\bs{B}\mid\bs{X}} \left( \bs{b} \mid \bs{x} \right)$, referred to as the likelihood,  is of the form
\begin{align*}
    \pi_{\bs{B}\mid \bs{X}} \left( \bs{b} \mid \bs{x} \right) \propto \exp \left( -\frac{1}{2} \left( \bs{b} - F \bs{x}  \right)^H (I_{N_T}\otimes \Gamma)^{-1} \left( \bs{b} - F \bs{x}  \right) \right).
\end{align*}
Without loss of generality, after applying a whitening transformation, we may assume that the noise term is zero mean white Gaussian so that  the likelihood becomes
\begin{align*}
    \pi_{\bs{B}\mid \bs{X}} \left( \bs{b} \mid \bs{x} \right) \propto \exp \left( -\frac{1}{2} \|  \bs{b} - F \bs{x}  \|_2^2 \right),
\end{align*} 
where for notational economy we do not introduce new symbols for the whitened quantities.
Let us assume that $\bs{x}$ admits a sparse representation in terms of a few columns of a matrix $W \in \mathbb{R}^{N \times P}$, referred to as a \emph{dictionary}, i.e, there is a vector $\bs{z}$ with only very few non-zero entries such that $\bs{X}=W \bs{Z}$, where $\bs{Z}$ is a $\mathbb{C}^P$-valued random variable. In that case we can write the likelihood in terms of the vector of coefficents $\bs{z}$ as
\begin{align} \label{eq:Likelihood}
    \pi_{\bs{B}\mid  \bs{Z}} \left( {\bs{b}} \mid \bs{z} \right) \propto \exp \left( -\frac{1}{2} \| {\bs{b}} - {F} W \bs{z} \|_2^2 \right).
\end{align}
In Section \ref{subsec:dictionary}, we outline how the dictionary matrix $W$ is built so that its columns capture the spatiotemporal characteristics of the unknown, and its rows judiciously represent the temporal behavior of each component.

\section{Derivation of the Spatiotemporal Dictionary}
\label{subsec:dictionary}

In a dictionary-based approach, the unknown of interest is represented as a linear combination of suitably chosen templates constituting the dictionary atoms. A preference is usually given to sparse linear combinations of the atoms, i.e., only a few atoms from a possibly overcomplete dictionary are needed to explain the unknown. In dictionary-based inverse problems the quantity to be estimated is the vector $\bs{z}$ of the weights of the dictionary atoms. To guarantee that only few atoms are needed to approximate the unknown, it is important that the atoms of the underlying dictionary capture well the salient features of the original unknown. In this subsection we outline how to form a dictionary $W$ for the dynamic inverse problem of interest so that the aforementioned goals are met.

We begin by considering the matrix $X$ to be estimated in two different ways, i.e., as the ensemble of its columns and as the ensemble of its rows, and we introduce the notation
\[
 X = \left[\begin{array}{cccc}
 \bs{x}^{(1)} & \bs{x}^{(2)} & \cdots & \bs{x}^{(N_T)}\end{array}
 \right] = 
 \left[
 \begin{array}{c} \bs{x}_{(1)}^T \\
 \bs{x}_{(2)}^T \\
 \vdots \\
 \bs{x}_{(N_p)}^T\end{array}
 \right] 
 \in{\mathbb R}^{N_p\times N_T}.
\]
Because of the organization of the entries of $X$, each  column vector $\bs{x}^{(t)}$ represents a snapshot of the unknown at time $t$, and the row vector $\bs{x}_{(j)}^T$ represents the time evolution of a specific spatially localized quantity, e.g., a pixel or voxel value. 

The design principle behind the construction of the dictionary $W$ that will be used to represent our unknown should reflect our a priori expectations. 

Consider first the spatial dimension. At each discrete time instance $t$, we assume that the snapshot $\bs{x}^{(t)}$ in the spatial direction  admits a representation in terms of an appropriately chosen spatial dictionary $S\in{\mathbb C}^{N_p\times N_S}$, thus 
\[
 \bs{x}^{(t)} = S \bs{\omega}^{(t)}, \quad 1\leq t\leq N_T,
\]
where $\bs{\omega}^{(t)}$ is a coefficient vector.
Exploiting the partitioning of $X$ into the ensemble of its columns, we can express the spatial prior in matrix-matrix form as
\begin{equation}\label{omega}
 X = S\Omega, \quad \Omega\in{\mathbb C}^{N_s\times N_T},
\end{equation}
where $\Omega$ is the coefficient matrix. The selection of the matrix $S$ is application-dependent, and will be discussed later.

Similarly, for a given spatial degree of freedom $j$, assume that the corresponding time trace $\bs{x}_{(j)}$ admits a representation in terms of a properly chosen dictionary $E$, i.e., we may write
\begin{equation}\label{E formula}
 \bs{x}_{(j)} = E\bs{\xi}_{(j)},\quad 1\leq j\leq N_s,
\end{equation}
for some coefficient vector $\bs{\xi}_{(j)}$.
Again, organizing the coefficient vectors as columns of a matrix $\Xi$, we write a matrix-matrix form of the dictionary ansatz as
\begin{equation}\label{xi}
    X^T = E\, \Xi, \mbox{ or } X = \Xi^T E^T
\end{equation}
where $\Xi$ is the corresponding coefficient matrix. To reconcile the two a priori models (\ref{omega}) and (\ref{xi}), we write an ansatz that is compatible with both of the assumptions,
\begin{equation}\label{final}
 X = S Z E^T,
\end{equation}
where $Z$ is assumed to be a matrix such that
\[
 Z E^T = \Omega, \quad S Z = \Xi.
\]
Notice that since the matrices $\Omega$ and $\Xi$ are not given a priori, the above identities will not be enforced. In fact, the matrices $\Omega$ and $\Xi$ are simply intermediate quantities of no interest in the algorithm.

Consider now the vectorized form of the representation (\ref{final}). A standard result about Kronecker products shows that the vectorized form  $\bs{x} = vec(X)$ of the unknown to be estimated is a linear transformation of the vectorized form $\bs{z}= vec(Z)$ of the matrix $Z$. Indeed, applying the operator $vec$ on both sides of (\ref{final}) yields
\begin{equation}\label{kron}
 \bs{x}  = \big(E\otimes S\big) \bs{z}.
\end{equation}
Since multiplying a vector by a matrix is equivalent  to taking a linear combination of the columns of the matrix scaled by the corresponding entries of the vector, (\ref{kron}) is a way to express the unknown to be estimated in terms of the columns of the matrix $E\otimes S$, suggesting that a judicious choice of the spatiotemporal dictionary in line with our a priori beliefs about the solutions is given by 
\[
 W = E\otimes S.
\]
Assuming furthermore that the dictionary $W$ is representative enough, we seek to express the vector $\bs{x}$ in terms of as few dictionary entries as possible, which means that we pursue a sparse representation.
In summary, we reformulate the dynamic inverse problem as the dictionary matching problem of finding a sparse representation for $\bs{x}$ of the form
\begin{equation}\label{z to x}
 \bs{x} = W\bs{z}, \mbox{ $\bs{z}$ sparse.}
\end{equation}

Finally, we remark that the dictionary $W$ does not need to be formed explicitly, because the proposed matrix-free algorithms will require only the action of the dictionary and its transpose on vectors, which can be carried out through the steps
\begin{equation}\label{matvec 1}
 \bs{z} \buildrel {vec^{-1}} \over \longrightarrow Z \longrightarrow SZE^T 
 \buildrel vec \over \longrightarrow  W\bs{z},
\end{equation}
and similarly, the action of $W^T$ can be carried out by observing that $W^T = E^T\otimes S^T$, so
\begin{equation}\label{matvec 2}
 \bs{y} \buildrel {vec^{-1}} \over \longrightarrow Y \longrightarrow S^T YE 
 \buildrel vec \over \longrightarrow W^T\bs{y}.
\end{equation}
Before discussing the problem-specific dictionaries, in the next section we review the Bayesian sparsity promoting algorithm that is the computational backbone of our proposed algorithm.

\section{Sparsity-Promoting Hierarchical Bayesian Model}
\label{sec:SPHB}
The Iterative Alternating Sequential Algorithm (IAS) is a sparsity-promoting  algorithm for the solution of linear inverse problems, based on a hierarchical  Bayesian framework. The IAS algorithm proceeds according to a block descent scheme, where at each iteration two minimization problems are solved, corresponding to a partitioning of the unknowns into two groups, formally with a similar logic to the scheme of the Alternating Direction Method of Multipliers (ADMM) \cite{boyd_distributed_2010}. Both the IAS and the ADMM algorithms require the setting of some hyperparameters, however in the IAS, very little hyperparameter tuning is required, owing to its derivation from the Bayesian framework. In this section we give a brief description of the IAS algorithm for computing a sparse approximation for the maximum a posteriori (MAP) estimate of the Bayesian inverse problem in the dictionary learning framework, and review some of its properties. A more in depth  description of the IAS algorithm in the context of dictionary matching can be found in \cite{pragliola2021overcomplete}. 

In the Bayesian framework, the unknown coefficient vector $\bs{z}$ in the representation (\ref{z to x}) is a random variable, denoted by $\bs{Z}$. Since we expect a priori that most of the entries of $\bs{Z}$ will vanish, with only a few at unknown locations different from zero, it is reasonable to assume that the entries of $\bs{Z}$ are mutually independent. We will assume that $\bs{Z}$ 
follows a Gaussian distribution with zero mean and diagonal covariance, thus   $Re (\bs{z}) \sim \mathcal{N} (0, {\rm diag}(\bs{\theta}))$ and $Im (\bs{z}) \sim \mathcal{N} (0, {\rm diag}(\bs{\theta}))$, where $\theta$ is the diagonal vector of the covariance matrix.   
Since the diagonal entries $\theta$, which are the variances of the components of $\bs{Z}$, are not known, they are modeled also as random variables, the random vector being denoted by $\bs{\Theta}$. Our prior assumes that the real and imaginary part of each component of $\bs{Z}$ have the same variance, $\Theta_i$, which, as will be shown later, amounts to promoting group sparsity of real-imaginary coordinate pairs. \\ 
The conditionally Gaussian prior for $\bs{Z}$ given $\bs{\Theta}=\bs{\theta} $ is of the form 
 \begin{eqnarray}
\pi_{\bs{Z} \mid \bs{\Theta}} (\bs{z} \mid \bs{\theta}) & \propto & \prod_{i=1}^{P} \theta_i^{-\frac{1}{2}} \exp{ \left( -\frac{1}{2}  \sum_{i=1}^{P} \left(\frac{Re\left( \bs{z}_i \right)^2 + Im \left( \bs{z}_i \right)^2}{\theta_i} \right) \right) } \nonumber\\
& \propto & \exp{ \left( -\frac{1}{2} \| L_{\theta}^{-1/2} \bs{z} \|_2^2 - \frac{1}{2} \sum_{i=1}^{P} \ln \left( \theta_i \right) \right) }\label{eq:CondGauss},
\end{eqnarray}
where $L_{\theta}$ is the Cholesky factor of the matrix ${\rm diag}(\bs{\theta})$. 
\\\\
In line with our prior beliefs, we need to define a distribution model for $\bs{\Theta}$, referred to as the \emph{hyperprior}, that ensures the positivity of its entries and naturally promotes sparsity in $\bs{z}$. One choice that is particularly advantageous from the point of view of computations is to assume that the entries of $\bs{\Theta}$ are mutually independent and identically distributed, and that each one follows the Gamma distribution
\begin{align*}
\pi_{\bs{\Theta}_i}(\theta_i) = \frac{\theta_i^{\beta - 1}}{\Gamma(\beta) \vartheta^\beta} \exp{\left( -\frac{\theta_i}{\vartheta} \right)}
\propto \exp{ \left( -\frac{\theta_i}{\vartheta}  + \left( \beta - 1 \right) \ln \left(  \theta_i \right) \right)},
\end{align*}
where $\Gamma$ is the Gamma function, and $\beta$ and $\vartheta$ are the shape and scale parameters respectively, which will both be constant for all components and are discussed further below. 
This choice of a fat tailed hyperprior ensures that all realizations of $\Theta_i$ are positive, with expectation close to zero, while allowing the possibility of a few large outliers. Entries of $\bs{z}$ with small variance will be penalized more heavily for deviating from their zero mean than entries with larger variance, thus promoting the sparsity of $\bs{z}$. In view of the independence of the entries of $\bs{\Theta}$ it follows that
\begin{eqnarray} \label{eqn:hyperprior}
\pi_{\bs{\Theta}}(\bs{\theta}) &\propto& \prod_{i=1}^{P} \exp{ \left( -\frac{\theta_i}{\vartheta}  + \left( \beta - 1 \right) \ln \left(  \theta_i \right) \right)} \nonumber \\
& \propto &  \exp{ \left( - \sum_{i=1}^{P} \left( \frac{\theta_i}{\vartheta}  + \left( 1 - \beta \right) \ln \left(  \theta_i \right) \right) \right)}. 
\end{eqnarray}
\\\\
The joint posterior distribution for $\bs{Z}$ and $\bs{\Theta}$, obtained from Bayes' formula and the law of total probability, is of the form :
\begin{eqnarray*}
\pi_{\bs{Z}, \bs{\Theta} \mid \bs{B}} \left( \bs{z}, \bs{\theta} \mid {\bs{b}} \right) &\propto& \pi_{\bs{B} \mid \bs{Z}, \bs{\Theta}} \left( {\bs{b}} \mid \bs{z}, \bs{\theta} \right) \pi_{\bs{Z}, \bs{\Theta}} \left( \bs{z}, \bs{\theta} \right) \\
& \propto &  \pi_{\bs{B} \mid \bs{Z}, \bs{\Theta}} \left( {\bs{b}} \mid \bs{z}, \bs{\theta} \right) \pi_{\bs{Z} \mid \bs{\Theta}} \left( \bs{z} \mid \bs{\theta} \right) \pi_{\bs{\Theta}} \left( \bs{\theta} \right).
\end{eqnarray*}
It follows from (\ref{eq:Likelihood}), (\ref{eq:CondGauss}) and (\ref{eqn:hyperprior}) that  
\begin{align*}
\pi_{\bs{Z}, \bs{\Theta} \mid \bs{B}} \left( \bs{z}, \bs{\theta} \mid {\bs{b}} \right) \propto \exp{\left( - f(\bs{z}, \bs{\theta}) \right)},
\end{align*}
where 
\begin{align} \label{eqn:Gibbs}
f(\bs{z}, \bs{\theta}) =  \frac{1}{2} \| {\bs{b}} - {F} W \bs{z} \|_2^2 + \frac{1}{2} \| L_{\bs{\theta}}^{-1/2} \bs{z} \|_2^2
+ \sum_{i=1}^{P} \left( \frac{\theta_i}{\vartheta}  + \left( \frac{3}{2} - \beta \right) \ln \left( \theta_i \right)) \right).
\end{align}
While the full Bayesian analysis of the posterior density would require sampling from this density, it is customary to summarize the posterior density by computing the MAP estimate, or the minimizer of the Gibbs energy $f(\bs{z}, \bs{\theta})$ given by (\ref{eqn:Gibbs}), 
\begin{align} \label{eqn:MAP}
(\bs{z}_{\rm MAP}, \bs{\theta}_{\rm MAP}) = \argmin_{\bs{z}, \bs{\theta}} \left( f(\bs{z}, \bs{\theta}) \right).
\end{align}
In the next section we review computational methods specifically designed for computing the MAP estimate for sparsity promoting hierarchical Bayesian models.

\section{The Iterative Alternating Sequential Algorithm}
\label{sec:IAS}

The IAS algorithm proceeds to minimize \eqref{eqn:Gibbs} according to a block descent scheme, where at each iteration two minimization problems are solved. A more in depth description of the IAS algorithm in the context of dictionary matching can be found in \cite{pragliola2021overcomplete}. 

The IAS algorithm is an iterative procedure for computing the minimizer of the functional (\ref{eqn:Gibbs}) by alternatively fixing $\bs{\theta}$ and minimizing with respect the $\bs{z}$, then fixing $\bs{z}$ to the updated value and minimizing with respect to $\bs{\theta}$. Each IAS iteration computes a new approximation of the optimizer in two steps, each solving an optimization problem with respect to one block of variables  while fixing the remaining ones at their most recent values. The IAS iterations continue until a stopping criterion is satisfied.  
\subsection{IAS stage one: Updating $\bs{z}$}
At the $k$th iteration round, assuming that the value of $\bs{\theta}= \bs{\theta}^{k-1}$ is given and fixed, we find the updated value of $\bs{z}$ by solving 
\begin{align} \label{eqn:IAS1}
\bs{z}^k = \argmin_{\bs{z}} \left( \frac{1}{2}\|\bs{b} - {F} W \bs{z} \|_2^2 + \frac{1}{2} \| L_{\theta^{k-1}}^{-1/2} \bs{z} \|_2^2 \right).
\end{align}
By writing the objective function in (\ref{eqn:IAS1}) in the form
\begin{align*}
\bs{z}^k = \argmin_{\bs{z}} \left( \frac{1}{2} \left\|  
    \begin{bmatrix} 
    {F} W \\
    L_{\bs{\theta}}^{-1/2}
    \end{bmatrix} \bs{z} - 
    \begin{bmatrix} 
    {\bs{b}} \\
    \bs{0}
    \end{bmatrix}
    \right\|_2^2
    \right),
\end{align*}
and introducing the new variable $\bs{\zeta} = L_{\bs{\theta}^{k-1}}^{-\frac{1}{2}} \bs{z}$, we can update $\bs{z}$ by solving 
\begin{align}
\bs{\zeta}^k = \argmin_{\bs{\zeta}} \left( \frac{1}{2} \left\|  
    \begin{bmatrix} 
    {F} W L_{\bs{\theta}}^{1/2}\\
    I
    \end{bmatrix} \bs{\zeta} - 
    \begin{bmatrix} 
    {\bs{b}} \\
    \bs{0}
    \end{bmatrix}
    \right\|_2^2
    \right), \; \; \bs{z}^k = L_{\bs{\theta}}^{1/2} \bs{\zeta}^k.
    \label{eg:LSksi}
\end{align}  
The solution of the minimization problem (\ref{eg:LSksi}) also solves the normal equation 
\begin{align} \label{eqn:Stage1Normal}
\left( A_{\theta}^H A_{\theta} + I \right) \bs{\zeta} = A_{\theta}^H {\bs{b}},
\end{align}
where $A_{\theta} = {F} W L_{\theta}^{1/2}$. 
In general, when the matrix $A$ is large,  it is not desirable to explicitly form the matrix $A_\theta$ or its inverse, but equation (\ref{eqn:Stage1Normal}) can be solved iteratively using  Krylov-subspace linear solvers, where the matrices $A_\theta$ and $A_\theta^H$ are needed only to compute matrix-vector products, which can be performed matrix-free, see formulas (\ref{matvec 1}) and (\ref{matvec 2}) and  further discussion in Section \ref{subsection:efficienctImplementation}. 
 
 \subsection{IAS stage two: Updating $\bs{\theta}$}
 After updating  $\bs{z}=\bs{z}^k$, phase 2 updates $\bs{\theta}$ to the minimizer of (\ref{eqn:Gibbs}) with respect to $\bs{\theta}$. Equivalently, we may write
\begin{align} \label{eqn:IAS2}
\bs{\theta}^k = \argmin_{\bs{\theta}} \left( \sum_{i=1}^{N_{\theta}} \left( \frac{\theta_i}{\vartheta}  + \eta \ln \left(  \theta_i \right) + \frac{(\bs{z}_i)^2}{2 \theta_i} \right) \right),
\end{align}
where, for convenience, we have defined $\eta = \frac{3}{2} - \beta$.
 
As it turns out, equation (\ref{eqn:IAS2}) admits an analytical solution. 
Because of the mutual independence of the components of $\bs{\theta}$, the minimizer of the summation can be computed component-wise by solving 
\begin{align*}
\theta_i^k = \argmin_{\bs{\theta}_i} \left( \frac{\theta_i}{\vartheta}  + \eta \ln \left(  \theta_i \right) -\frac{z_i^2}{2} \theta_i^{-1} \right).
\end{align*}
It follows from the first order optimality condition that the minimizer $\theta_i$ must satisfy
\begin{align*}
\frac{1}{\vartheta} - \eta \theta_i^{-1} - \frac{z_i^2}{2\theta_i^2} = 0.
\end{align*}
Multiplying both sides of the equation by $\theta_i^2$, yields a quadratic equation which has the unique nonnegative solution
\begin{align} \label{eqn:Stage2Analytic}
\theta_i^k = \frac{\vartheta}{2} \left( \eta + \sqrt{\eta^2 + 2 \frac{z_i^2}{\vartheta}} \right), 
\end{align}
thus completing the second phase of the updating process.

A pseudocode for solving the optimization problem (\ref{eqn:MAP})  by the IAS algorithm is given by Algorithm \ref{alg:IAS}. The first step is solved using the LSMR method. LSMR requires forward and adjoint operations which are shown in Algorithms \ref{alg:ForwardOp} and \ref{alg:AdjointOp} respectively.

\begin{algorithm}
\caption{IAS}\label{alg:IAS}
\begin{algorithmic}
\State $\bs{\zeta}_{(0)} = \left[  \; 0 \; \hdots \; 0 \right]^T$
\State $\bs{\theta}_{(1)} = \left[ \vartheta \; \vartheta \; \hdots \; \vartheta \right]^T $
\For{$k=1:max\_iter$}
    \State Step 1 : Solve $\left( A_{\bs{\theta}_{(k)}}^H A_{\bs{\theta}_{(k)}} + I \right) \bs{\zeta}_{(k)} = A_{\bs{\theta}_{(k)}}^H \hat{\bs{b}}$ using LSMR
    \State Step 2 : Compute $\bs{\theta}_{(k+1)} = \vartheta \left( \frac{\eta}{2} + \sqrt{\frac{\eta^2}{4} + \frac{\bs{\zeta}_{(k+1)}^2}{2 \vartheta}} \right)$
\If{$\frac{\| \bs{\theta}_{(k)} - \bs{\theta}_{(k-1)} \|}{\| \bs{\theta}_{(k-1)}\| } < tol$}
\State{break}
\EndIf
\EndFor
\end{algorithmic}
\end{algorithm}

\begin{algorithm}
\caption{Forward Operation}\label{alg:ForwardOp}
Matrix-free computation of $\bs{y} \rightarrow A_{\theta} \bs{y}$, where $A_{\bs{\theta}} = F \left( E \otimes S \right) L_{\bs{\theta}}^{\frac{1}{2}}$, and the constituent matrices are never explicitly formed.
\begin{algorithmic}
\State $\bs{y} \rightarrow \sqrt{\bs{\theta}}  \circ \bs{y}$ \; \; Component-wise multiplication
\State $Y = \text{reshape} (\bs{y}, N_S, N_E$)
\State $Y \left( : , i \right) \rightarrow S ( Y \left( : , i \right) )$ \; \; Apply $S$ to each column of $Y$ (in parallel)
\State $Y \left( i , : \right) \rightarrow E (Y \left( i , : \right) )$  \; \;  Apply $E$ to each row of $Y$ (in parallel)
\State $\bs{y} = vec(Y)$
\State $\bs{y} \rightarrow \hat{F} \bs{y}$ Apply forward operator
\end{algorithmic}
\end{algorithm}
\begin{algorithm}
\caption{Adjoint Operation}\label{alg:AdjointOp}
 Matrix-free computation of $\bs{y} \rightarrow A_{\theta}^H \bs{y}$, where $A_{\bs{\theta}}^H =  L_{\bs{\theta}}^{\frac{1}{2}} \left( E^H \otimes S^H \right) F^H$, and the constituent matrices are never explicitly formed.
\begin{algorithmic}
\State $\bs{y} \rightarrow \hat{F}^H \bs{y}$ \; \; Apply adjoint forward operator
\State $Y=\text{reshape} \left( \bs{y}, N_P, N_T \right)$ \; \; Prepare for Kronecker product
\State $Y[:, i] = S^H \left( Y[:, i] \right)$ \; \; Apply $S^H$ to each column of $Y$ (in parallel) 
\State $Y[i, :] = E^H \left( Y[i, :] \right)$ \; \; Apply $E^H$ to each row of $Y$ (in parallel) 
\State $\bs{y} = vec \left( {Y} \right)$
\State $\bs{y} \rightarrow \bs{y} \circ \sqrt{\bs{\theta}}$
\end{algorithmic}
\end{algorithm}

\subsection{Hyperparameters}
The hypermodel for the variance parameter $\bs{\theta}$ contains two parameters, $\beta$ and $\vartheta$ that need to be selected by the user. The advantage of the IAS algorithm is that these hyperparameters have a clear interpretation, making their selection intuitive, and furthermore, the algorithm is not overly sensitive to the choices. \\
As discussed in \cite{calvetti2023bayesian}, the parameter $\eta = \beta - 3/2>0$ controls the sparsity of the solution, and it has been shown that at the limit $\eta\to 0+$, the IAS solution approaches the $\ell^1$-penalized solution. In general, the algorithm is not very sensitive to the choice of $\eta$ as long as the parameter is small. This observation is corroborated by the numerical sensitivity analysis discussed later in this work. 
Recalling that the expectation of the Gamma distribution is $\beta\vartheta \approx 3\vartheta/2$ for small $\eta$,  we see that the second hyperparameter $\vartheta$ is related to the expected value of the variances $\theta_i$, and therefore a rule of thumb is to choose $\vartheta$ small enough to guarantee that most variances $\theta_i$ are small. Again, the algorithm is not sensitive to the specific choice, as will be demonstrated by the numerical experiments.

\subsection{Implementation Considerations} 
\label{subsection:efficienctImplementation}
By far the most computational time in each IAS iteration is spent in stage 1, that is, in solving the normal equations (\ref{eqn:Stage1Normal}). Special attention needs to be paid to the implementation to increase the computational efficiency when the problem is of high dimensionality. The Least Squares Minimal Residual (LSMR) algorithm is the method of choice for the solution of large scale linear and least squares problems. As pointed out earlier, the application of LSMR requires that the matrix-vector products with the dictionary and its adjoint are implemented efficiently.

\section{Results}
\label{sec:results}
In this section, we demonstrate the proposed dynamic reconstruction approach on two real-world dynamic imaging experiments, using Computed Tomography (CT) and Dynamic Contrast-Enhanced MRI (DCE-MRI) respectively. 
In these applications, the design of the dictionary $W$ outlined in the previous section is based on specific prior assumptions concerning the spatio-temporal behavior of the unknown as stated below.

\begin{itemize}
\item {\bf{A priori assumptions about the spatial behavior.}} 
 At each discrete time instance $t$, we assume that the snapshot $\bs{x}^{(t)}$ in the spatial direction  admits a sparse representation in terms of an appropriately chosen spatial dictionary $S\in{\mathbb C}^{N_p\times N_S}$.  In both computed examples, our choice of spatial dictionary matrix $S$ is the inverse 2D Haar wavelet transform with symmetric boundary conditions and three scales, implemented using the PyWavelets Python package \cite{lee_pywavelets_2019}. 
\item {\bf{A priori assumptions about the time traces.}}
The expectation that most space locations will remain mostly constant over time, with a few possibly large changes at a few time instances, can be expressed by representing each time evolution vector $\bs{x}_{(j)}$ as a sparse sum of increments, i.e., 
\begin{equation}\label{L formula}
 L \bs{x}_{(j)} = \bs{\xi}_{(j)}, \quad 1\leq j\leq N_p,
\end{equation}
where $L\in{\mathbb R}^{N_T\times N_T}$ is the invertible first order finite difference matrix,
\[
 L = \left[\begin{array}{rrrr}
  1 & & & \\
 -1 & 1 & & \\
    &\ddots & \ddots & \\
    &       & -1     & 1\end{array}\right],
\]    

Comparing formulas (\ref{E formula}) and (\ref{L formula}),
we see that a consistent choice for the matrix $E$ is
\[
 E = L^{-1} = \left[\begin{array}{llll}
 1 & & & \\ 1 & 1 & & \\ 
 \vdots & & \ddots & \\
 1 & 1 & \cdots & 1\end{array}
 \right],
\]
Observe that it is possible to compute the action of the matrix $E$ on a vector without the need to form it explicitly through a cumulative summation.  For alternative ways of forming spatio-temporal dictionaries, see, e.g., \cite{bubba_regularization_2024}).
\end{itemize} 

We solve the test problems by using the dictionary-based IAS algorithm, and for comparison, using a standard ADMM sparse solver, described in detail in the Appendix.

The following stopping tolerances were all set to $10^{-8}$: $tol$ for IAS as used in Algorithm $\ref{alg:IAS}$, $\epsilon_{rel}$ and $\epsilon_{abs}$ for ADMM as they are defined in \cite{boyd_distributed_2010}, and $ATOL$ and $BTOL$ for LSMR as they are defined in \cite{fong_lsmr_2011}. Stopping tolerances were also applied to maximum iteration counts for IAS, ADMM and the LSMR solvers, and their values are specified for each numerical experiment.
\subsection{Computed Tomography}
\begin{figure}[h]
        \centering
        \includegraphics[width=\textwidth]{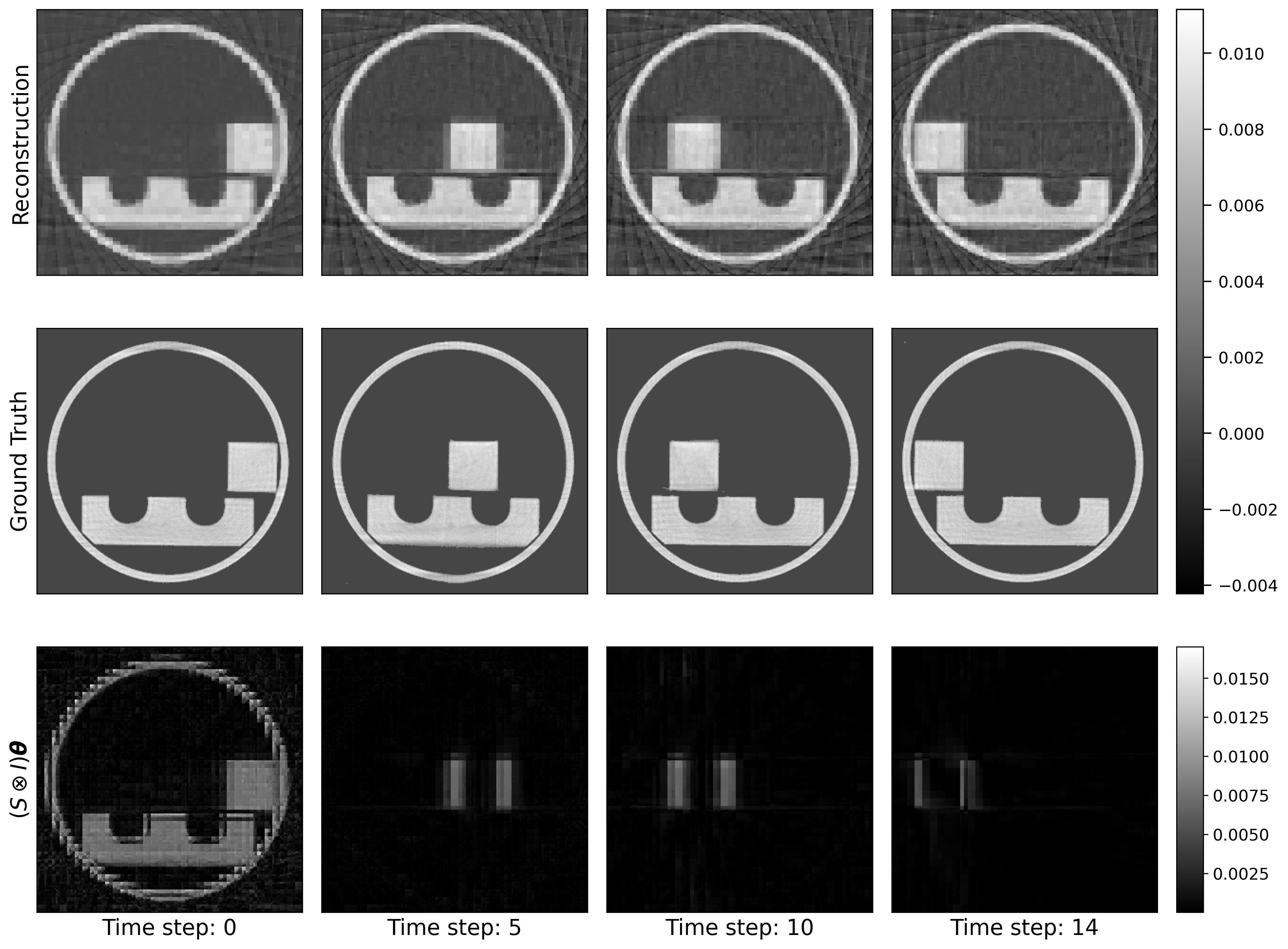}
        \caption{CT reconstructions with optimised hyper-parameter selection (top), compared to the ground truth (middle). The variable, $\theta$, has been transformed to image space via application of the inverse wavelet transform at each timestep.}
        \label{fig:CTIASOptRecon}
\end{figure}

The CT data originates from the open source STEMPO dataset \cite{heikkila_stempo_2023}, specifically the file "stempo\_seq8x180\_2d\_b8.mat". We give a brief description of the experimental set-up and data; see \cite{heikkila_stempo_2024} for full details. The data are collected from a micro-CT device with an X-ray source and detector on opposite sides of the target. The target is placed on a platform which rotates between acquisitions, so that projections at multiple angles are measured. The target consists of a short piece of $80\,$mm diameter polypropene pipe containing a static object made of high density polyethylene (HDPE). A dynamic HDPE object is also placed inside the pipe and moved across the pipe diameter in a straight line, under the control of a step motor during data acquisition (see ground truth in Figure \ref{fig:CTIASOptRecon}). Note that this is an ideal target for compressed sensing reconstructions, as a perfect reconstruction is possible with sparse increments.
The projection data are acquired at 2 degree intervals around the target, making for a total of 90 measurements per 180$^\circ$ rotation. The data are partitioned into 16 parcels of full 180$^\circ$ rotations representing the snapshots, for a total of $N_T=16$ time points. To generate a strongly undersampled test case, each of the 16 parcels is downsampled by selecting every $5^{th}$ projection angle, resulting in $18$ evenly spaced projection angles per timepoint, and a total of 288 projections over all timepoints. The resulting alternating sampling pattern repeats every two timepoints, where odd timepoints are sampled in the range $[0, 180)$ and even timepoints in the range $[180, 360)$. A ground-truth comparison is created from the full dataset (90 projection angles per timepoint) by solving a lightly Tikhonov-regularized ($\lambda = 10^{-3}$) least-squares problem. The ground truth is thresholded by removing values below $2.9 \times 10^{-3}$, and pixels outside the pipe cross section are set to zero. Note that the same masking is performed when computing structural similarity scores in the analysis to avoid biasing the metric. The forward operator $F^{(t)}$ at each time point $t$  is the X-ray projection operator for a 2D fan beam geometry and was implemented using the ASTRA toolbox version 2.4.0 \cite{van_aarle_fast_2016}. IAS and ADMM reconstructions are performed on the CT dataset using a maximum of 10 outer iterations and a maximum of 50 inner (LSMR) iterations, with the exception of the results shown in Figure \ref{fig:CTRuntimes}.

\begin{figure}[h]
        \centering
        \includegraphics[width=\textwidth]{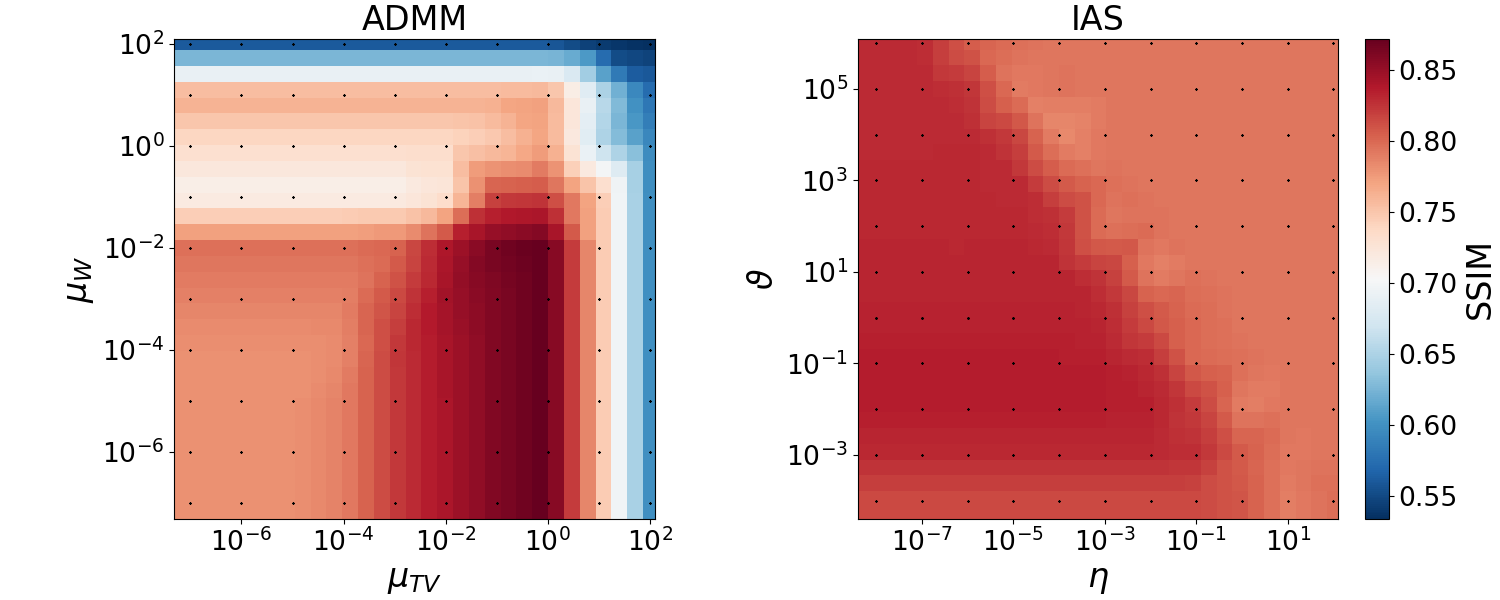}
        \caption{Heatmaps depicting the SSIM scores for ADMM (left) and IAS (right) reconstructions in terms of the selected hyper-parameter pairs. The black dots indicate the location of measured values over which linear interpolation was performed.}
        \label{fig:CTGriddedSSIM}
\end{figure}

\begin{figure}[h]
        \centering
        \includegraphics[width=\textwidth]{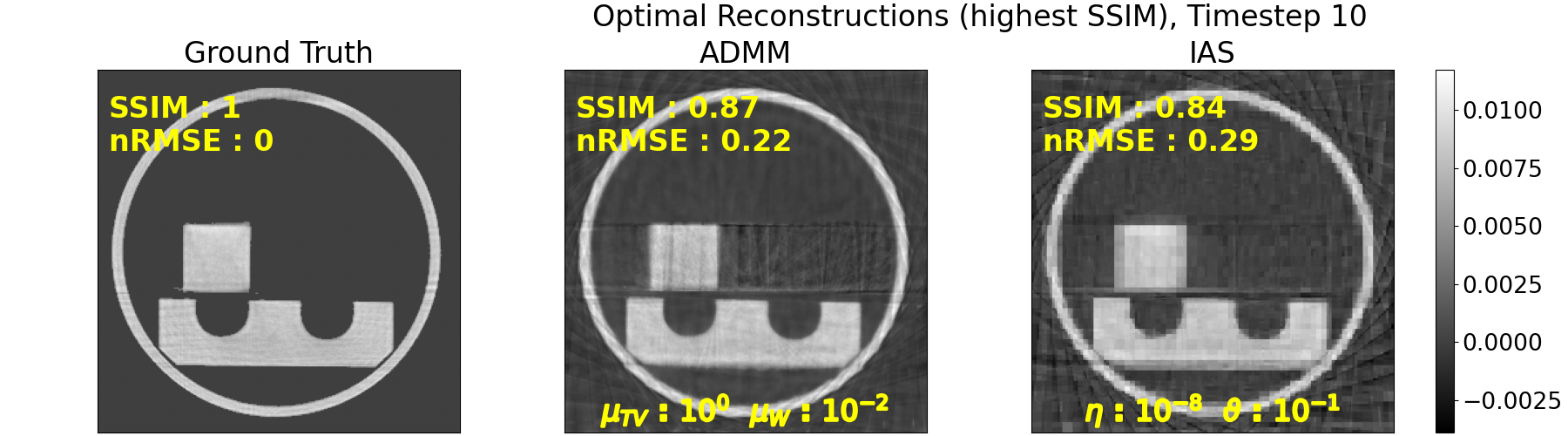}
        \includegraphics[width=\textwidth]{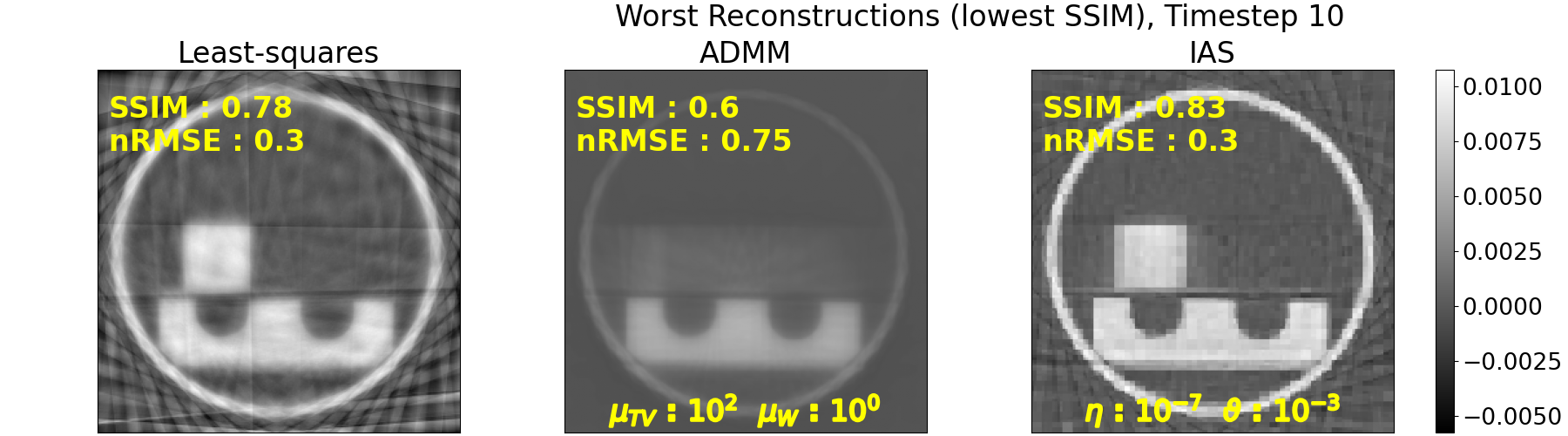}
        \caption{Optimal (top) and worst-case (bottom) reconstructions for ADMM and IAS, respectively. Worst-case reconstructions are chosen as those with the lowest SSIM among hyper-parameter settings within two orders of magnitude of the optimal values. The ground truth and least-squares reconstructions are shown on the left for reference, and SSIM and normalised root-mean-squared error (nRMSE) scores for each reconstruction are displayed in yellow.}
        \label{fig:SSIMWorstBestRecons}
\end{figure}

\begin{figure}[h]
        \centering
        \includegraphics[width=\textwidth]{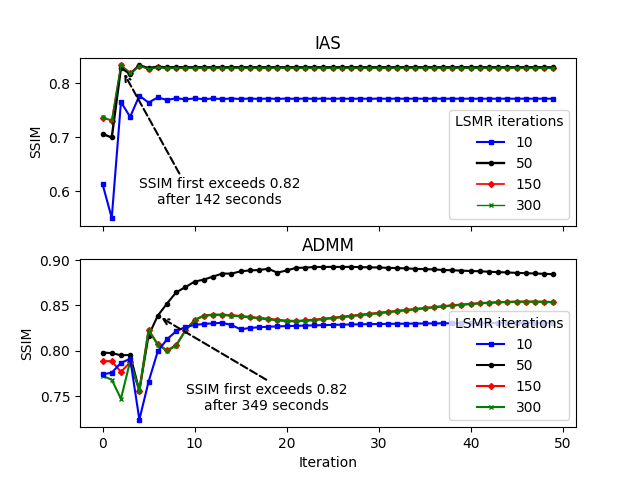}
        \caption{SSIM scores for the IAS (top) and ADMM (bottom) reconstructions at each iteration of the respective algorithms. Plots are shown using a stopping criteria of 10, 50, 150, and 300 LSMR iterations respectively. Note that the black, red and green lines overlap in the IAS plot; and the red and green lines overlap in the ADMM plot.}
        \label{fig:CTRuntimes}
\end{figure}

The first row of Figure \ref{fig:CTIASOptRecon} shows four time points of an IAS reconstruction, performed using the optimized hyper-parameter pair, $\eta=10^{-8}$, $\vartheta=10^{-1}$, which was selected based on the results shown in Figure $\ref{fig:CTGriddedSSIM}$. There is a good matching between the reconstructions and the ground truth shown in the second row. A slight blurring on the left and right edges of the moving block, i.e., along the direction of motion can be observed, likely due to the choice of temporal prior. The third row shows the inverse wavelet transform of the estimated $\bs{\theta}$, i.e., $\left( S \otimes I \right) \bs{\theta}$, at four time points. Recalling that $\bs{\theta}$ represents a vector of variances associated with the dictionary coefficients $\bs{z}$, we expect it to match the sparsity pattern of $\bs{z}$, and for the sparsity pattern of $\left( S \otimes I \right) \bs{\theta}$ to match the sparsity pattern of $\left( S \otimes I \right) \bs{z}$. At each of the time points 5, 10, and 14 we see two vertical lines, corresponding to the locations of the left and right edges of the moving block where the time derivative of the image is high, and hence where we expect $\left( S \otimes I \right)\bs{z}$ to be large.  
The time point 0 on row 3 demonstrates the effect of our implicit temporal boundary condition $\bs{x}_t = 0$ for $t<0$, meaning effectively that at timestep $1$, the increment conicides with the image itself, and hence we expect the first time point of $\left( S \otimes I \right)\bs{z}$ to resemble the reconstructed image. \\\\
Figure \ref{fig:CTGriddedSSIM} shows the sensitivity of the reconstructed image quality to the choice of hyper-parameter pairs for the IAS and ADMM respectively. Image quality was assessed using the Structural Similarity Index (SSIM) \cite{zhou_wang_image_2004}, which was time-averaged over all images in the reconstructed time series. The SSIM was computed using the {\em metric.structural\_similarity} function from the scikit-image toolbox \cite{scikit_image} and, due to the large number of zeros in the ground truth image, was tuned to favor contrast and structural similarities over luminance using the input parameters $K_1=0.1$, $K_2=0.1$, $sigma=1.5$ and $win\_size=11$, see the cited reference for description of the parameters. In both panels, the SSIM scores are plotted for model-specific hyper-parameter values all varying over eight orders of magnitude. The SSIM scores for IAS are all within the range 0.78-0.84 with a large flat region of SSIM scores near 0.84, while the ADMM SSIM scores are in the range 0.54-0.87 and there is a sharp drop off from the optimum. The SSIM scores for IAS are highest when $\eta$ is low, and this behavior is expected given the expected sparsity of the dictionary coefficient vector for this target. While a slightly larger SSIM score can be achieved with ADMM with careful tuning and with the knowledge of the ground truth, in practice when no ground truth is available, it is difficult or impossible to find the required hyper-parameters to achieve this score. In contrast, the insensitivity of the IAS to hyper-parameter pairs, allows the user to deviate far from the optimum and still achieve a high quality reconstruction, reflecting robustness of the IAS for the hyperparameter selection.\\\\
Figure \ref{fig:SSIMWorstBestRecons} shows the IAS and ADMM reconstructions that achieve the maximum SSIM values in Figure \ref{fig:CTGriddedSSIM} (the optimal reconstructions), as well as worst-case reconstructions defined as those with the lowest SSIM scores among hyper-parameter pairs within two orders of magnitude of the corresponding optimal values. The unregularized least-squares reconstruction and ground truth are shown on the left for reference. While the optimal IAS and ADMM reconstructions appear similar, the worst ADMM reconstruction is of significantly poorer quality than the worst IAS one, with the moving block seen as a low intensity blur across the entire path of its motion. This again highlights the robustness of IAS to hyper-parameter selection in contrast to CS ADMM. \\\\
Figure \ref{fig:CTRuntimes} demonstrates how the SSIM score changes over each iteration of the IAS (top) and ADMM (bottom), respectively, using a varied stopping tolerance on the number of LSMR iterations for the linear solve step. The 10, 50, and 150 iteration cases all reached this tolerance, while the 300 case generally terminated around 220 iterations due to the other stopping criteria in place. Note that, for both algorithms, the linear solve takes up the majority of computational time, and hence the run times are largely determined by the accuracy of the linear solve (number of LSMR iterations) and the number of outer iterations. In all four cases, the IAS reaches a steady SSIM score after 6-8 iterations. The SSIM scores are almost identical after three IAS iterations between the 50, 150 and 300 LSMR iteration cases, and systematically slightly lower when only 10 LSMR iterations are used. This finding indicates robust convergence with respect to the number of LSMR iterations, as going beyond a certain threshold in the LSMR iterations does not essentially affect the results. 

For the ADMM reconstructions, the SSIM scores are higher for 50 LSMR iterations than for 10 or 150 or 300. The unexpected fact that the increase in accuracy by increasing the LSMR iterations beyond 50 results in a decrease in reconstruction quality in terms of the similarity score could suggest that the truncation of the linear solver iterations at 50 introduces an inherent regularization favorable for the reconstruction quality, making the selection of the LSMR iterations a more delicate issue in the ADMM than it is in IAS.

Finally, we point out that with ADMM, only the 10 iteration solution reaches a steady SSIM score within 50 iterations of the algorithm, while for larger iteration indices, the curves continue to grow before reaching a saturation level. Moreover, in the optimal case with 50 LSMR iterations, the best performance is observed around 25 iterations, followed by a descending trend in the quality. This observation is a further indication that the selection of the number of iterations for ADMM is a more delicate issue than for IAS. Furthermore, while ADMM reaches a slightly higher peak SSIM score, it requires much longer run times than IAS, as well as a prior knowledge on how to truncate the iterations on the linear solver. For comparison, in both subplots, an arrow points to the first iteration for the 50 LSMR solver where the SSIM score surpasses 0.82. At this point, IAS run time is less than half of that of ADMM.
\subsection{Dynamic Contrast-Enhanced MRI (DCE-MRI)}
The in-vivo DCE-MRI dataset was collected with a 9.4T pre-clinical Agilent MRI device from an adult female Wistar rat with a C6 glioma, with Gadovist contrast agent injected into the blood stream through the tail vein over 3 seconds, beginning one minute into the acquisition. The acquisition employed a gradient echo based radial golden-angle sampling scheme with 128 points per radial spoke. The total number of radial spokes in the time series data was 7310, corresponding to approximately 5 minutes scan time. For full details on the dataset, see \cite{DCEAnatomicalPrescan}.

For the dynamic reconstruction, the data was split to 215 frames of 34 spokes per image, corresponding to approximately 1.3 seconds per frame with approximately $12\times$ sub-Nyquist sampling at the periphery of $k$-space. A time series of 215 $128\times 128$ images were reconstructed using the IAS, ADMM and unregularized least squares methods. The forward model is a block diagonal Non-Uniform Fast Fourier Transform (NUFFT) operator, implemented using the {\em jax-finufft} library from the Flat Iron Institute \cite{barnett_aliasing_2021, barnett_parallel_2019}. Both the IAS and ADMM reconstructions were performed using a stopping criteria of 300 inner (LSMR) iterations and 50 outer iterations. This example allowed for more iterations than the CT example owing to the larger dimensionality of the data.

Figure 5 displays IAS and ADMM reconstructions, each using two different hyper-parameter pairs, with two order of magnitude shift in $\vartheta$ and $\mu_T$, respectively. Reconstructions are shown at time points 30 and 100, before and after injection of the Gadovist agent. In all reconstructions, at time point 100, one can see increased pixel intensity in the upper middle and right side of the image where the tumor is located. The IAS reconstructions with different hyper-parameter pairs appear similar to each other, and similar to the left-hand side ADMM reconstruction. However the right-hand side ADMM reconstruction, which used a smaller temporal regularisation strength, appears blurred and noisy.
\begin{figure}[h]
        \centering
        \includegraphics[width=\textwidth]{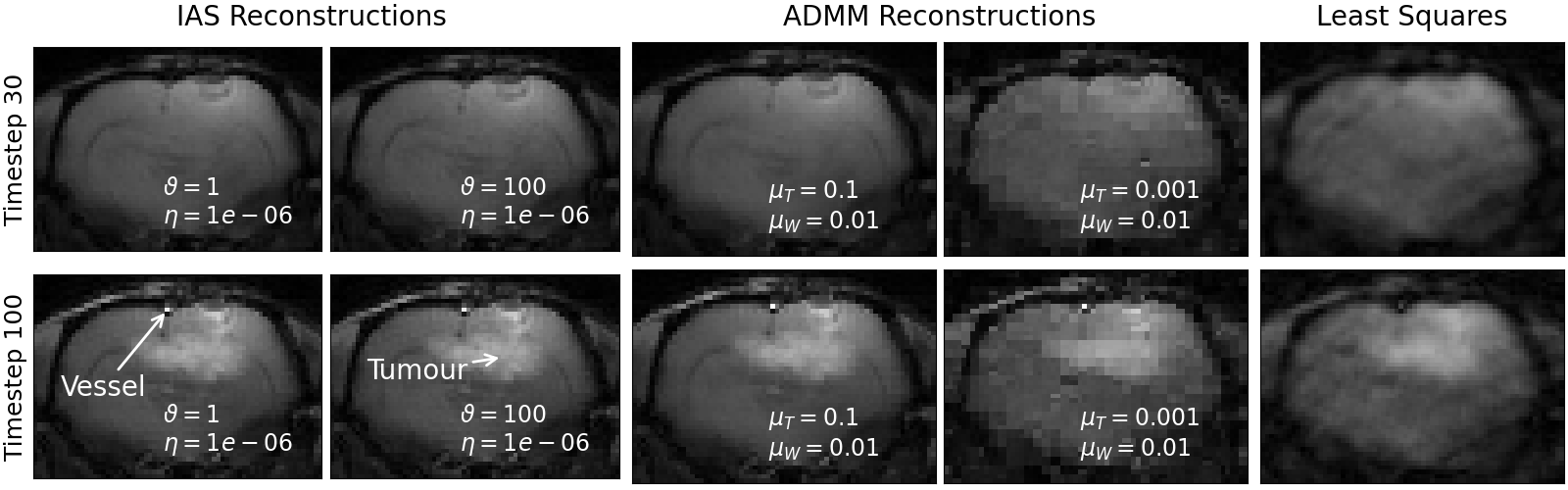}
        \caption{Reconstructed DCE-MRI images at timepoints 30 (top) and 100 (bottom). IAS (left) and ADMM (middle) reconstructions are shown for two different hyper-parameter pairs, which are indicated in white on the images. White arrows on the IAS images indicate the locations of pixels in the vessel and tumour regions. In the absence of a ground truth reference, the unregularized least squares reconstruction is shown on the far right.}
        \label{fig:DCEReconstructions}
\end{figure}

Figure 6 shows single pixel signal time courses, located in a vessel and in the tumor, respectively. The white arrows in Figure 5 point to the locations of each respective pixel. Both the IAS and ADMM tumor time courses exhibit a dip due to the first pass of the bolus, followed by rising pixel intensity due to accumulation of the contrast agent to the glioma through the broken blood-brain barrier. The time courses with the different IAS hyper-parameter pairs are almost identical, while the ADMM with the same relative change (two orders of magnitude) in temporal regularization parameter leads to notably different solutions, exhibiting large noise (right, orange) with the smaller regularization. The IAS reconstructions appear smooth, while the more heavily regularised ADMM reconstruction (right, blue) is jagged, possibly due to the staircasing effect. The IAS time courses in the vessel pixel are also similar to each other, with the green line (lower $\vartheta$ value) reaching a slightly lower amplitudes. The ADMM vessel time courses also show possible staircasing in the line with higher temporal regularization, and heavy noise in the case of smaller regularization.
\begin{figure}[h]
        \centering
        \includegraphics[width=\textwidth]
        {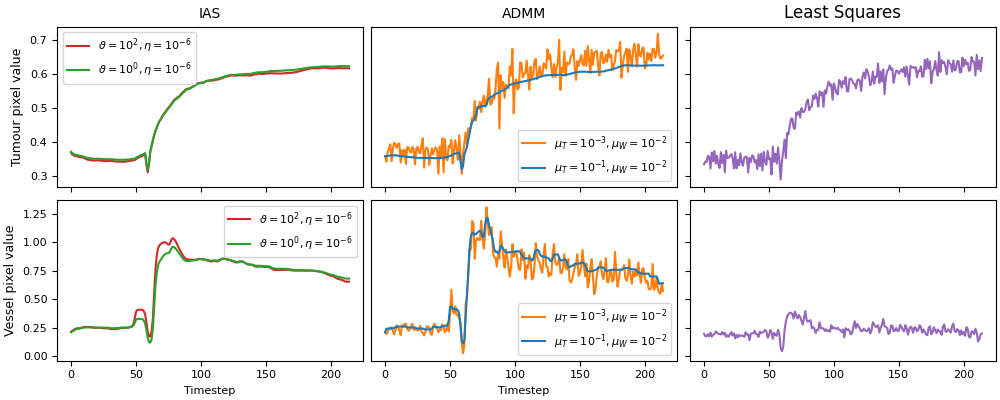}
        \caption{Single pixel time courses in the pixels indicated by the white arrows in Figure \ref{fig:DCEReconstructions}. The tumour time courses (top) are shown for the IAS (left) and ADMM (middle) with two different hyper-parameter pairs each. The vessel time courses are displayed on the bottom two plots. The unregularized least squares reconstruction is shown on the far right.}
        \label{fig:DCESignals}
\end{figure}

\section{Discussion and Conclusions}
\label{sec:conclusions}
We introduced a novel approach to dynamic inverse problem reconstruction using a spatio-temporal dictionary-based  sparsity-promoting hierarchical Bayesian model. The viability of the proposed approach was demonstrated on two real-world datasets, and the results were compared to those obtained with a commonly used compressed sensing formulation with a standard ADMM-based solver. 
The two algorithms produced reconstructions of similar quality, however, the IAS reconstructions had much lower sensitivity to hyper-parameter selection. While in the presence of the ground-truth in the X-ray tomography example, it was possible to tune the hyperparameters of the ADMM-based solution to reach a slightly higher similarity index than the proposed method, it is important to note that normally no ground-truth is available as in the DCE-MRI experiment, making the insensitivity of the IAS to hyperparameter selection an important asset. Moreover, the X-ray tomography target is ideal for the $\ell^1$-sparsity reconstruction underlying the CS algorithm, while in the MRI example, it was clear that non-optimal selections of the hyperparameters in ADMM lead either to noisy reconstructions (under-regularization) or to staircase artifacts (over-regularization) while the IAS solution is void of these artifacts.

Due to the similar structure of both algorithms, the difference in run times comes down to the number of algorithmic (outer) iterations and the number of (inner) iterations used by the linear solver. In the CT case, the IAS algorithm requires far fewer outer iterations than ADMM to reach a steady reconstruction accuracy. The IAS reconstruction quality was also unchanged when the maximum number of the inner LSMR iterations was increased beyond 50, while the ADMM reconstruction quality behaved inconsistently, leading to decreased reconstruction quality when the inner iteration were increased beyond 50 iterations. The phenomenon is similar to the semi-convergence of iterative solvers for ill-posed problems, and it is typically hard to detect when no ground truth is available.

The IAS was formulated under the assumption of Gaussian measurement noise, which is not necessarily true for MRI data \cite{gudbjartsson_rician_1995}. 
While the Gaussian approximation worked well with the real data considered in this article, a careful analysis of the noise distribution can be incorporated in the IAS, leading typically to a non-quadratic likelihood model, see, e.g., \cite{bardsley_hierarchical_2010} for an example of non-Gaussian noise model. 
\appendix

\section*{Appendix: Baseline Comparison : The Alternating Direction Method of Multipliers}
\label{sec:ADMM}
As a baseline comparison, we present a Compressed Sensing (CS) formulation for solving a dynamic inverse problem with the forward model (\ref{eqn:ForwardModel}). The compressed sensing problem is formulated as follows:
\begin{align} \label{eqn:CSProblem}
\bs{x^*} = \argmin_{x} \left( c\left( \bs{x} \right) \right),
\end{align}
where the cost function $c \left( \bs{x} \right)$ is:
\begin{align*}
c \left( \bs{x} \right) = \frac{1}{2} \| \hat{\bs{b}} - \hat{F} \bs{x} \|_2^2 + \mu_1 \| H_1 \bs{x} \|_1 + \mu_2 \| H_2 \bs{x} \|_1,
\end{align*}
where the $H_1$ and $H_2$ are transforms such that $H_1 \bs{x}_1$ and $H_2 \bs{x}_2$ are expected to be sparse and the functions $h(\bs{x}_1) = \| H_1 \bs{x}_1 \|_1$ and $g(\bs{x}_2) = \| H_2 \bs{x}_2 \|_1$ are convex. The positive constants $\mu_1$ and $\mu_2$ control the weightings of the penalty terms relative to the data misfit, and require tuning to the given problem. 
Equation (\ref{eqn:CSProblem}) can be reformulated as the following constrained minimisation problem:
\begin{align*}
\bs{x^*}, \bs{q^*} = \argmin_{\bs{x}, \bs{q}} \left(  \frac{1}{2} \| \hat{\bs{b}} - \hat{F} \bs{x} \|_2^2 +\| \bs{q} \|_1 \right),
\end{align*}
such that,
\begin{align*}
H \bs{x} - \bs{q} = 0,
\end{align*}
where $H = \begin{bmatrix} \mu_1 H_1 \\ \mu_2 H_2 \end{bmatrix}$. To solve this constrained minimisation problem with ADMM, we first form the augmented Lagrangian:
\begin{align*}
\mathcal{L}_{\rho} \left( \bs{x}, \bs{q}, \bs{u} \right) = \frac{1}{2} \| \hat{\bs{b}} - \hat{F} \bs{x} \|_2^2 + \| \bs{q} \|_1 + \frac{\rho}{2} \| H \bs{x} - \bs{q} + \bs{u} \|_2^2 - \frac{\rho}{2} \| \bs{u}\|_2^2,
\end{align*}
where $\bs{u}$ is the so-called scaled dual variable.
Now Equation (\ref{eqn:CSProblem}) can be solved by the following iterative update procedure, known as the ADMM:
\begin{align*}
\bs{x}^{(k+1)} &= \argmin_{\bs{x}} \left(  \mathcal{L}_{\rho} \left( \bs{x}, \bs{q}^{(k)}, \bs{u}^{(k)} \right) \right) = \argmin_{\bs{x}} \left( \frac{1}{2} \| \hat{\bs{b}} - \hat{F} \bs{x} \|_2^2 + \frac{\rho}{2} \|H \bs{x} - \bs{q}^{(k)} + \bs{u^{(k)}} \|_2^2 \right), \\
\bs{q}^{(k+1)} &= \argmin_{\bs{q}} \left(  \mathcal{L}_{\rho} \left( \bs{x}^{(k+1)}, \bs{q}, \bs{u}^{(k)} \right) \right) = \argmin_{\bs{q}} \left( \| \bs{q} \|_1 + \frac{\rho}{2} \| H \bs{x}^{(k+1)} - \bs{q} + \bs{u^{(k)}} \|_2^2 \right), \\
\bs{u}^{(k+1)} &= \bs{u}^{(k)} + \rho \left( H \bs{x}^{(k+1)} - \bs{q}^{(k+1)} \right).
\end{align*}
Let $S_{\kappa}$ be the soft thresholding operator. The full algorithm can be stated as:
\begin{algorithm}
\caption{ADMM}\label{alg:ADMM}
\begin{algorithmic}
\While{not converged}
    \State Solve $\left( F^T F + \rho H^T H \right) \bs{x}^{(k+1)} = F^T \hat{\bs{b}} + \rho H^T \left( \bs{q}^{(k)} - \bs{u}^{(k)} \right)$ for $\bs{x}^{(k+1)}$
    \State Compute $\bs{q}^{(k+1)} = S_{\frac{1}{\rho}} \left( H \bs{x}^{(k+1)} + \bs{u}^{(k)} \right)$
    \State Compute $\bs{u}^{(k+1)} = \bs{u}^{(k)} + \rho \left( H \bs{x}^{(k+1)} - \bs{q}^{(k+1)} \right)$
\EndWhile
\end{algorithmic}
\end{algorithm}

\section*{Acknowledgments}
This project has received funding from the European Union’s Horizon 2020 research and innovation programme under the Marie Skłodowska-Curie agreement No101034307. The work of VK was supported by Research Council of Finland (RCF) Flagship of Advanced Mathematics for Sensing, Imaging and Modeling grant (no.358944), the RCF Centre of Excellence in Inverse Modelling and Imaging grant (no.~353084). The work of DC was also partly supported by the NSF grant DMS 2513481 and by a Collaboration Grant by the Simons Foundation, and that of ES by the NSF grants DMS 2204618 and DMS 2513481.
All computations were performed on servers provided by UEF Bioinformatics Center, University of Eastern Finland, Finland.
 \bibliographystyle{siam} 
\bibliography{references}
\end{document}